\documentclass[11,reqno]{amsart}

\usepackage{amssymb,latexsym}
\usepackage{amsmath,amsthm}
\usepackage{enumerate}
\usepackage[mathscr]{eucal}
\usepackage{graphics} 

\def\phi{\varphi}
\def \P {\mathbb P}
\newtheorem{Th}{Theorem}
\begin{document}

\centerline{\bf TESTS BASED ON CHARACTERIZATIONS,}
\centerline{\bf AND THEIR EFFICIENCIES: A SURVEY}

\bigskip

\centerline{ Nikitin Ya. Yu.}

\bigskip

\noindent{Abstract.} {\footnotesize A survey of goodness-of-fit and symmetry tests based on the characterization properties of distributions is presented. This approach became popular in recent years. In most cases the test statistics are functionals of $U$-empirical processes. The limiting distributions and large deviations  of new statistics under the null hypothesis are described. Their local Bahadur efficiency for various parametric alternatives is calculated and compared with each other as well as with diverse previously known tests. We  also describe new directions of possible research in this domain.}

\section{Introduction}

This survey is dedicated to the statistical tests based on characterizations. This is a relatively new idea which manifests growing popularity in the context of goodness-of-fit and symmetry testing. The idea to build goodness-of-fit tests using the cha\-racterizations of distributions belongs to Yu. V. Linnik \cite{Lin1}. At the end of this wide-ranging paper he wrote: \emph{ "... one can raise the issue of the construction of goodness-of-fit tests for testing composite hypotheses based on the equal distribution of the two relevant statistics $ g_1 (x_1, ... x_r) $ and
$ g_2 (x_1, ... x_r), $ and on the reduction of such question to the homogeneity tests."} This sentence was the guiding star which showed the researchers the right direction in the new and unexplored domain.

Currently, in the world literature there exist hundreds of various characterizations of probability distributions, see, e.g., \cite{Kagan}, \cite{Galamb}, \cite{Kotz}, and \cite{KKM}. Many characterizations according to Linnik's idea imply the corresponding statistical tests. Such tests are attractive because they reflect some \emph{intrinsic and hidden properties} of probability distributions connected with the given characterization, and therefore can be more efficient or more robust than others.

Moreover, one should keep in mind that any hypothesis has to be tested with several possible criteria. The point of the matter is that with absolute
confidence we can only reject it, while each new test which fails to reject the null-hypothesis gradually brings the statistician closer to the perception that this hypothesis is true. We find it pertinent to quote here the famous assertion by Einstein \cite{Einst}:\emph{ "No amount of experimentation can ever prove me right; a single experiment can prove me wrong."} Hence, we are interested in building new statistical tests based on novel ideas, specifically using the characterizations.

But the theory of such tests is intricate, and the study of their asymptotic pro\-perties including limiting behavior, and especially their asymptotic efficiency began only after 1990. Before that there existed few exceptions like the paper \cite{Va}, of which later Mudholkar and Tian \cite{Mu} wrote: \emph{"Va\v{s}icek (1976) was the first to recognize that the characterization results can be logical starting points for developing goodness-of-fit tests."}

\smallskip

\hrule width 190pt height 0.4pt depth 0.4pt

\medskip

{\footnotesize
\noindent \emph{  2010 Mathematics Subject Classification} 60F10,\ 62F03,\ 62G20,\ 62G30.\\
\noindent \emph{ Key words and phrases}. Characterization of distribution, goodness-of-fit, symmetry, $U$-statistics, Bahadur efficiency.
}
\newpage

Probably these authors were unfamiliar with the seminal paper by Linnik cited above who was surely the first to propose the idea under discussion.
In the abstract of the paper \cite{Hash} published in 1993, Hashimoto and Shirahata proposed one of the first tests of fit based on characterizations and wrote: \emph{"However, since no test statistics based on characterizations are known, our test will be worth considered."} This citation shows that in the beginning of 1990-s the tests based on characterizations were unusual and sparse. But since that time the state of affairs changed significantly. Numerous new tests based on characterizations were build, and their study \emph{gradually  acquired the traits of a theory.} We want to trace an outline of this theory and its main achievements within the last 25 years.

We begin by general constructions explaining the structure of tests used in this domain. Next we pass to concrete problems like testing of exponentiality, normality or symmetry, and describe the main developments of last period of time. We are \emph{mainly interested in the asymptotic efficiency of our tests} though the results of power simulation are also possible and interesting. At the end of the paper we pose some problems and trace new directions of research. In most cases, the proofs are omitted, otherwise this survey would exceed the size of the paper in a journal.

\section{$U$-statistics and $U$-empirical distributions}

Let $X_1, X_2,\ldots,X_n$  be i.i.d. observations with continuous df $F$. We begin by testing the \emph{composite goodness-of-fit} problem
$$H_{0}: F \in  \mathfrak{F}, $$ where ${\mathfrak F}  $ is some family of df's, against the alternative $$H_{1}: F\notin { \mathfrak F}.$$
Typical examples are testing  exponentiality, normality or symmetry of a sample.

Next exposition will be based on $U$-statistics and their variants. Currently $U$-statistics play an important role in Statistics and Probability. They appeared in the middle of 1940-s in problems of unbiased estimation, but after the crucial paper of Hoeffding \cite{Hoeff} it became clear that the numerous valuable statistics are just $U$-statistics (or von Mises functionals having very similar asymptotic theory.) The most complete exposition of this theory can be found in the monographs \cite{KorBor} and \cite{Lee}.

We consider $U-$statistics of the form
$$
    U_n={n \choose m}^{-1}\sum_{1\leqslant i_1<\ldots<i_m\leqslant
    n}{\Psi(X_{i_1},\ldots,X_{i_m})},\qquad n\geqslant m ,
$$
where $X_1, X_2,\dots$ is a sequence of i.i.d. rv's with common df $F$, while the kernel $\Psi: {\mathbb R}^m \to {\mathbb R}^1$ is a measurable symmetric function of $m \ge 1$ variables. The number $m$ is called the \emph{degree} of the kernel. We assume that the kernel $\Psi$ is integrable on $R^m$ and denote
$$
\theta(F) := \int...\int_{{\mathbb R}^m} \Psi(x_1,\ldots, x_m) dF(x_1)... \ dF(x_m).
$$
In the sequel we need the notations
$$
    \psi(x):= \mathbb{E}_F\{\Psi(X_1,\ldots,X_m)|X_1 = x\}, \quad
    \Delta^2 := \mathbb{E}_F\psi^2(X_1) - (\theta(F))^2.
    $$

The function  $\psi$ is called the one-dimensional \emph{ projection} of the kernel $\Psi$ and plays an important role in asymptotic theory. If $\Delta^2 > 0$ that specifies the so-called non-degenerate case, the limiting distribution of $U-$statistics is normal as discovered by Hoeffding \cite{Hoeff}. He proved that
if $ \mathbb{E}_F \Psi^2(X_1,\ldots,X_m) < \infty$ and $\Delta^2 > 0,$ then as $n \to \infty$ one has convergence in distribution
 \begin{equation}
 \label{normal}
    {\frac{\sqrt{n}}{m\Delta}}\left(U_n - \theta(F)\right)\mathrel{\stackrel{\makebox[0pt]{\mbox{\normalfont\tiny d}}}{\longrightarrow}} N(0,1).
\end{equation}

Consider, in conformity with Linnik, the characterization of some probability law by the identical distribution  of two statistics $g_1(X_1, \ldots , X_r)$ and $g_2(X_1, \ldots , X_s).$ The examples of such characterizations will be given further.

We can build  two \emph{ $U$-empirical df's}
\begin{gather*}
L^1_n(t)={n \choose r}^{-1}\sum_{1 \leq i_1<\ldots <i_r \leq n}\textbf{1}\{g_1(X_{i_1}, \ldots , X_{i_r})< t\}, \quad t\in {\mathbb R}^1,\\
L^2_n(t)={n \choose s}^{-1}\sum_{1 \leq i_1<\ldots <i_s \leq
n}\textbf{1}\{g_2(X_{i_1}, \ldots , X_{i_s})<t\},\quad t\in {\mathbb R}^1.
\end{gather*}

The theory of $U$-empirical df's was developed in 80-s, see, e.g. \cite{helmers}, \cite{Jan}, \cite{KorBor} and is similar to the theory of usual empirical df's. By Glivenko-Cantelli theorem for {$U$-empirical df's} we have (wp 1) as $n\to \infty:$
$$
\begin{array}{ll}
\cr L^1_n(t) \rightrightarrows L^1(t):=P(g_1(X_1, \ldots , X_r)<t), \\
\cr L^2_n(t)\rightrightarrows L^2(t):=P(g_2(X_1, \ldots , X_s)<t).
\end{array}
$$

Under $H_0 $  for large $n$  we have $ L^1_n(t) \approx L^2_n(t)$ so that we can use this closeness for goodness-of-fit testing.
Over much of this survey we  consider two types of statistics: the integral one
$$\label{I_n}
 I_n=\int_{{\mathbb R}^1} \left(L^1_n(t)-L^2_n(t)\right)dF_n(t),
$$
where $F_n(t) $ is the usual empirical df, and of Kolmogorov type, namely
$$\label{D_n}
 D_n= \sup_{t \in {\mathbb R}^1 }\mid L^1_n(t)-L^2_n(t)\mid.
$$

Such statistics can have rather different behavior depending on the type of cha\-racterization and underlying distribution, accordingly the statistical tests based on them can have distinct limiting properties, power and efficiency.

\section{Outline of Bahadur theory}

Suppose that we want compare two sequences of statistics $I_n$ and $D_n$ by their asymptotic efficiency. Among many types of efficiencies, see \cite[Ch.1]{Niki}, we select the Bahadur efficiency because, unlike Pitman efficiency, it can be calculated for statistics with non-normal limiting distribution. This is the primary reason to use it in the present context as the Kolmogorov type statistics
have non-normal limiting distributions. Hodges-Lehmann efficiency has other drawbacks, in particular, it does not discriminate two-sided tests,
see, e.g., \cite[Ch.1]{Niki}. In this section we shortly describe main points of Bahadur theory, see the complete exposition in \cite{Bahad} and \cite{Bah2}.

Let $s=(X_1, X_2,\ldots)$ be a sequence of i.i.d. rv's with the distribution $ P_\theta, \theta \in \Theta,$ on $(\mathcal{X},\mathcal{A})$. We are testing the null-hypothesis $$ { H_0: \theta \in \Theta_0 \subset \Theta \subset {\mathbb R}^1}$$
 against the alternative $${ H_1: \theta \in \Theta_1 =\Theta \setminus \Theta_0.}$$

 For this problem we use the sequence of test statistics
$T_n(s)=T_n(X_1,\ldots, X_n).$ The Bahadur approach  prescribes one to fix the power of concurrent tests and to
compare the exponential rates of decrease of their sizes for  the
increasing number of observations and fixed alternative. This
exponential rate for a sequence of statistics $\{T_n\}$ is usually
proportional to some non-random function $c_T(\theta)$ depending
on the alternative parameter $\theta \in \Theta_1$ which is called the \emph{
exact slope} of the sequence $\{T_n\}$.
The Bahadur asymptotic relative efficiency (ARE) $\,
e_{V,T}^{\, B} (\theta)$ of two sequences of statistics
$\,\{V_n\}$ and $\,\{T_n\}$ is defined by means of the formula
$$e_{V,T}^{\,B}(\theta) = c_V(\theta)\,\big/\,c_T(\theta)\,.$$

The exact slope can be found by the  fundamental theorem of Bahadur \cite{Bahad}:
\begin{Th}
\emph{ Suppose that the following two
conditions hold:
\[
\hspace*{-3.5cm} \mbox{a)}\qquad  T_n \
\stackrel{\mbox{\scriptsize $P_\theta$}}{\longrightarrow} \
b(\theta),\qquad \theta > 0,\nonumber \] where $-\infty <
b(\theta) < \infty$, and $\stackrel{\mbox{\scriptsize
$P_\theta$}}{\longrightarrow}$ denotes convergence in probability
under $P_{\theta}$.
\[
\hspace*{-2cm} \mbox{b)} \qquad \lim_{n\to\infty} n^{-1}  \ln \P_{H_0} \left( T_n \ge t \right) \ = \ - h(t)\nonumber
\] for any $t$ in an open interval $I,$ on which $h$ is
continuous and $\{b(\theta), \: \theta > 0\}\subset I$. Then
$c_T(\theta) \ = \ 2 \ h(b(\theta)).$}
\end{Th}

Often the exact Bahadur ARE is uncomputable for any alternative depending on
$\theta$ but it is possible to calculate the local Bahadur ARE as
$\theta \in \Theta_1$ approaches the null-hypothesis. Then one speaks about the
\emph{ local} Bahadur efficiency and \emph{ local} Bahadur exact slopes \cite{Niki}.

Let $K(\theta, \theta_0) \equiv K(P_\theta,P_{\theta_0})$ be the Kullback-Leibler distance between $P_\theta$ and $P_{\theta_0},$ see, e.g.,
\cite{Bah2} or \cite{VDV}. Put for any $\theta \in \Theta_1$
\begin{eqnarray*}
K(\theta,\Theta_0):=\inf\{K(\theta,\theta_0):\theta_0 \in \Theta_0\}.
\end{eqnarray*}

The Bahadur - Raghavachari inequality (the analog of Cram\'{e}r - Rao inequality
in testing), see \cite{Bahad}, \cite{Niki}  states that for any $\theta \in \Theta_1$ one has
\begin{eqnarray*}
  c_{T}(\theta)\leq 2K(\theta,\Theta_0).
\end{eqnarray*}

Hence we may define the (absolute) \emph{ local Bahadur efficiency}
of the sequence  $\{T_n\}$ by the formula
\begin{eqnarray*}\label{eff}
eff_T = \lim_{\theta\to \partial
\Theta_0}{c_{T}(\theta)/2K(\theta,\Theta_0)}.
\end{eqnarray*}

Only in exceptional cases
$$
c_{T}(\theta)= 2K(\theta,\Theta_0), \quad \forall \theta \in \Theta_1.
$$

Therefore one can be interested in those $F= {\mathbb {L}}(X_1)$ for which
\begin{eqnarray*}
{eff_T = \lim_{\theta\to \partial
\Theta_0}{c_{T}(\theta)/2K(\theta,\Theta_0)}=1.}
\end{eqnarray*}
We  call this property the \emph{ local optimality} in Bahadur sense.
 An interesting question is to describe those alternatives for which the considered tests are locally optimal in Bahadur sense. The idea ascends to Bahadur \cite{Bah2}  but was developed by the author, see \cite{NN}, \cite[Ch.6]{Niki} and subsequent papers, e.g., \cite{NiPe}. However, we leave this direction apart as it requires considerable space to enounce the obtained results.

First condition of Theorem 1 is a variant of the Law of Large Numbers and its verification is easy. On the contrary, the second condition of this theorem describes the rough (logarithmic) large deviation asymptotics of test statistics under the null-hypothesis and is non-trivial. To verify it, we often use the theorem on large deviations of $U$-statistics by Nikitin and Ponikarov \cite{Nikiponi}:
\begin{Th}
Let $V_n$ be a sequence of $U$-statistics with centered, bounded and  non-degenerate kernel $\Psi.$ Then
\begin{eqnarray}
\label{NIP}
\lim_{n \to \infty} n^{-1}\ln \P \{V_n\geq a\}=-\sum_{j=2}^{\infty}b_j a^j,
\end{eqnarray}
where the series with numerical coefficients $b_j$ converges for sufficiently small $a > 0$, and $b_2=(2m^2 \Delta^2)^{-1},$ where $\Delta^2$
is the variance of the projection of the kernel $\Psi.$
\end{Th}

Large deviations for the {supremum} of the family of non-degenerate $U$-statistics  $\sup_{t \in T} U_n(t),$ where $U_n(t)$ for each $t \in T$  is a $U$-statistic with the non-degenerate kernel $\Xi(X,Y;t)$  which corresponds to Kolmogorov type statistics,  were studied in \cite{nikiNS}. The result is similar to (\ref{NIP}) but slightly more involved.

\section{Desu's characterization and corresponding tests of exponentiality}

One of most simple characterizations of exponential distribution belongs to Desu, see \cite{Desu}:
\begin{Th}
 \emph{ Let $X$ and $Y$ be non-negative i.i.d. rv's with df  differentiable at zero. Then
$X \stackrel{d} = 2\min (X,Y)$ iff $X$ and $Y$ are exponentially distributed.}
\end{Th}

  Using this characterization we will show how to build and analyze the corresponding tests of exponentiality.

Let $X_1,\dots,X_n$  be i.i.d. observations with non-degenerate df $F,$ and  let  $F_n$ be the corresponding empirical df.  We are testing the composite hypothesis \emph{
$$
 H_{0}: F(x)  \, \mbox{ is the df of exponential law with the density }\ f(x)=\lambda e^{-\lambda x}, x \geq 0,
$$
}
where $\lambda>0$ is some unknown parameter, against the alternative $H_1$ under which the hypothesis $H_0$ is wrong.

In this case we need the $U$-statistical empirical df $H_n$ which is defined as
$$
H_n(t)={n \choose 2}^{-1}\sum_{1 \leq i < j \leq n}\textbf{1}\{2\min
(X_{i},X_{j})<t\},\ t \geq 0.
$$
We will study two statistics  $$ I_n = \int_{0}^{\infty}(F_n(t)- H_n(t)) dF_n(t),$$ and
$$
D_n = \sup_{t \geq 0} |F_n(t)-H_n(t) |.
$$
Clearly their distribution under the null-hypothesis does not depend on $\lambda.$

The statistic  $I_n$  is asymptotically equivalent to the $U$-statistic of degree $3$ with the centered kernel
\begin{align*}
&\Psi(X,Y,Z)= \frac12 -\frac{1}{3}[\textbf{1}\{2\min
(X,Y)<Z\}-\\&-\textbf{1}\{2\min
(Y,Z)<X\}-\textbf{1}\{2\min
(X,Z)<Y\}].
\end{align*}

The projection of this kernel is $$E[\Psi(X,Y,Z)|Z=t] :=
\psi( s)=\frac13 e^{-s}-\frac{1}{18}-\frac{4}{9}e^{-3s},
$$
and the variance of the projection is
$\Delta^2:=E\psi^2(Z)= \frac{11}{3780}\approx 0.003.$

  By Hoeffding's theorem, see \cite{Hoeff} we get
\begin{Th}
Under the hypothesis $H_0$ one has convergence in distribution
$$
\sqrt{n}I_n \stackrel{d}{\longrightarrow} \mathcal{N}(0, 9\Delta^2), \quad \text{as} \quad   n \rightarrow \infty .
$$
\end{Th}

As to the large deviations, in our case we get for $a>0$
$$
\lim_{n\to \infty} n^{-1} \ln \P (I_n >a) = - f_I (a),
$$
where the function $f_I$ is continuous for sufficiently small $a>0,$ and, moreover,
$$
f_I (a) =  \frac{210}{11} a^2( 1 + o(1)), \, \mbox{as} \, \, a \to 0.
$$

By way of an example let calculate the local Bahadur efficiency of $I_n$ for the  Weibull alternative.
This means that the alternative df of observations is
$$ F(x,\theta) = 1- \exp(-x^{1+\theta}), x \ge 0, \theta \ge 0.$$

We find after some simple calculations that  as $\theta \to 0$
\begin{gather*}
c_I(\theta) \sim b_I(\theta)^2/(9\Delta^2) \sim 1.147 \theta^2.
\end{gather*}
The Kullback-Leibler  distance  \, $K(\theta)$ between $H_0$ and $H_1$   satisfies
$$
K(\theta)  \sim \pi^2\theta^2/12,  \, \theta \to
0.
$$
The \emph{local Bahadur ARE } of our test is consequently equal to
$$
eff(I):=\lim_{\theta \to 0}\frac{c_I(\theta)}{2K(\theta)}\approx {0.697}.
$$

Consider now the Kolmogorov-type statistic $D_n.$ The difference $F_n(t) -
H_n(t)$ is a \emph{family} of $U$-statistics with the kernels depending on $t\geq 0:$
\begin{equation}
\label{kern}
\Xi(X,Y;t)=
\frac12(\textbf{1}\{X<t\}+\textbf{1}\{Y<t\})
-\textbf{1}\{2\min (X,Y)<t\}).
\end{equation}

The limiting distribution of the sequence  $D_n$ is unknown. Critical values for statistics  $D_n$
can be found via simulation.

In our case the family  $\{\Xi(X,Y;t), t\geq 0\}$ from (\ref{kern}) is centered, bounded, non-degenerate and hence satisfies all conditions of Theorem 2.4 from \cite{nikiNS} on large deviations of
$U$-empirical Kolmogorov statistics. Therefore as $a>0$ by \cite{nikiNS}
$$
\lim_{n\to \infty} n^{-1} \ln P (D_n >a) = - f_{D} (a),
$$
where
$$
f_{D} (a) = 2a^2(1 + o(1)), \,\mbox{as} \, \, a \to 0.
$$

Consider again the   Weibull alternative.
Arguments similar to the case of integral statistic, see \cite{nikiNS}, show that
the  local Bahadur efficiency of the sequence $D_n$ is equal to  ${0.158}.$
We see that this efficiency is low and considerably smaller than in the integral case. It is a rule that Kolmogorov-Smirnov type statistics are less efficient than integral ones. There exist some exceptions but they are rare.

\section{Tests of exponentiality based on characterizations}

There are numerous characterizations of exponential law, probably more than of any other probability law, see, e.g., \cite{AH}, \cite{Az}, \cite{Bala} and \cite{Galamb}. We consider only few typical examples where the tests of fit are build and studied.

\subsection{Lack of memory property and corresponding tests}

First, we mention the celebrated "lack of memory" property which consists in that only the exponential distribution satisfies the functional equation in df's
$$
1- F(x+y) - (1-F(x) ) (1- F(y)) =0  \quad   \forall \, x, y \ge 0.
$$
Replacing $F$ by empirical df $F_n$, one obtains some empirical field, and the functionals of it can be used as test statistics for exponentiality,
see as examples of many papers in this direction \cite{AA}, \cite{Koul} and \cite{HKh}.

The "lack of memory" property can be simplified. Denote, following Angus \cite{Ang}, the class ${\bf D}_1$ of
right-continuous df's $F$ with $F(0-)=0$ and
$$
\lim_{h\to 0}\frac{F(h)-F(0)}h =l\in [0,\infty].
$$
Let $\bar F(x) = 1-F(x)$. Angus used the following statement
that belongs to Arnold and Gupta: the functional equation
$$
\bar F(2x)=\bar F^2(x) \quad \forall \ x \ge 0
$$
characterizes the exponential distribution in the class of such distributions in ${\bf D}_1$ which are not concentrated at 0. He
introduced a Kolmogorov type test based on this characterization and studied its properties in \cite{Ang}. Later its local Bahadur efficiency against standard alternatives was calculated in \cite{Nik} and \cite{nikiNS}. It turned out to be rather low.

\subsection{Characterizations based on order statistics}

   Another example is given by  Riedel-Rossberg characterization  in terms of order statistics, see \cite{RoRie}. Denote, as usually,
$X_{k,n}$ the $k$-th order statistic in the sample of size $n$, $1 \le k \le n.$ Then the following characterization holds.
\begin{Th}
{Two statistics  $X_{2,3}-X_{1,3}$ and $\min (X_1,X_2)$ are identically distributed iff the sample
$X_1, X_2, X_3$  consists of exponential rv's.}
\end{Th}

The construction of tests based on this characterization and their asymptotic analysis is performed similarly to the case of Desu characterization, see \cite{Vol}.

\medskip

Next consider the Ahsanullah's characterization. Suppose that the df $F$ belongs to the class of df's $\mathbb F_1$, where the failure rate function $f(t)/(1-F(t)) $ is monotone for $t \geq 0.$ Ahsanullah  \cite{Ahh} proved some characterizations of exponentiality within the class  $\mathbb F_1.$ We consider here only one of his  characterizations.
\begin{Th}
\emph{ Let $X$  and $Y$  be non-negative i.i.d. rv's from class $\mathbb F_1.$ Then  ${
|X-Y| \stackrel{d} =  2\min(X,Y)} $ iff $X$ and  $Y$ are exponentially distributed. }
\end{Th}

Corresponding tests were build and analyzed by Nikitin and Volkova in \cite{NikiVolki}.

\subsection{Characterization of Arnold and Villase\~{n}or}

Recently Arnold and Viilase\~{n}or  \cite{Arno} expressed in the form of hypothesis the following characterization of exponentiality:

\emph{Let  $X_1, X_{2}, \ldots $ be non-negative i.i.d. rv's with the density $f$ having derivatives of all orders around zero.  Then for any $k \geq 2$
$$
 \max(X_1, X_2, \ldots, X_k) \stackrel{d} = \sum_{i=1}^{k}\frac{X_i}{i}
 $$
iff $f$ is exponential.}

Arnold and Villase\~{n}or were able to prove this hypothesis only for $k=2$. Later  Yanev and Chakraborty \cite{Yan} proved it is true for $k=3$
 and later in \cite{Yan2} proved it for arbitrary $k,$ see also \cite{MilO}. Tests of exponentiality based on these characterizations and their efficiencies were studied in \cite{Jov} and in \cite{Volkk}.

 Other tests based on characterizations of exponential distribution in terms of order statistics were build and studied in \cite{NV}
and \cite{Milo}. One can mention also the characterization of exponential law by the same distribution of $X$ and $|X-Y|$ where $X,Y$ are i.i.d. rv's having absolute continuous distribution, see \cite{PR}. Some steps toward using it for testing were made in \cite{nikiNS}.

\subsection{Table of efficiencies}

Now we present a table of local Bahadur efficiencies of the majority of tests of exponentiality described above. We will compare them with well-known classic scale-free tests of exponentiality based on Greenwood statistic $R_n$, Moran statistic $M_n$ and Gini statistic $G_n.$ We recall that
$$
R_n =  2-\frac{1}{n}\sum_{i=1}^{n}\left(\frac{X_i}{\overline{X}}\right)^2,\, M_n = \frac{1}{n}\sum_{i=1}^n
\ln\left(\frac{X_i}{\overline{X}}\right)+\mathbf{C},\, G_n = \frac{\sum_{i,j=1}^{n}|X_i-X_j|}{2n(n-1)\overline{X}},
$$
where $\mathbf {C}$ denotes the Euler constant. We consider also the famous Lilliefors statistic \cite{Lil} which has the form
$$
Li_n=\sup_{x\geq 0} \big\vert 1-F_n(x)-e^{{-x}/{\overline{X}}} \big\vert,
$$
and belongs to Kolmogorov type statistics with estimated parameters.
On efficiencies of these statistics see \cite{NT}, \cite{NTchi}, \cite{TT}.

\medskip

  We consider the following standard alternatives against exponentiality:
\vspace{0.5cm}

i) { Weibull alternative} with the density
$$(1+\theta)x^\theta \exp(-x^{1+\theta}),\theta \geq 0, x\geq 0;$$

ii)  { Makeham alternative} with the density
$$(1+\theta(1-e^{-x}))\exp(-x-\theta( e^{-x}-1+x)),\theta \geq 0, x\geq 0;$$

iii) { linear failure rate alternative} with the density
$$(1+\theta x)e^{-x-\frac{1}{2}\theta x^2}, \theta \geq 0, x\geq 0;$$

Now let compare the values of local Bahadur efficiency for various statistics. All them are collected in Table 1 below and were calculated according to the approach developed above for the tests based on Desu characterization. The superscripts \emph{Ross} and \emph{Ahs} denote the statistics based on Riedel-Rossberg's or Ahsanullah's characterization.

 \begin{table}[!hhh]\centering
 \caption{Local efficiencies of tests for exponentiality.}
\begin{tabular}{|c|c|c|c|}
\hline
Statistic & Alternative  & Alternative &  Alternative  \\
& Weibull &  Makeham & linear failure rate \\
\hline
\multicolumn{4}{|c|}{Integral type statistics} \\
\hline
$I_n^{Ross}$ & 0.650 & 0.450 & 0.119 \\
\hline
$I_n^{Ahs}$ & 0.795 & 0.692 & 0.257 \\
\hline
 Gini & 0.876 & 1 & 0.750\\
\hline
Moran & 0.943 & 0.694 & 0.388  \\
\hline
Greenwood & 0.608 & 0.750 & 1  \\
\hline
\multicolumn{4}{|c|}{Kolmogorov type statistics} \\
\hline
$D_n^{Ross}$ & 0.320 & 0.207 & 0.047  \\
\hline
$D_n^{Ahs}$ & 0.450 & 0.470 & 0.187\\
\hline
Angus & 0.158 & 0.187& 0.073 \\
\hline
Lilliefors & 0.538 & 0.607 & 0.356 \\
\hline
\end{tabular}
\end{table}

We see that our tests based on characterizations are competitive with respect to other tests of exponentiality, the more that the alternatives were taken almost at random. However, the Gini test reaffirms its high reputation.

\section{Tests of normality}

Characterizations of normality are also numerous and mathematically content-rich. They are described in \cite{Ahsann}, \cite{Kagan}, \cite{Mathai}, and \cite{Bryc}, apart from many articles. We discuss here only few papers based on selected characterizations.

\subsection{Polya characterization}
One of first characterizations in the history of Statistics belongs to Polya \cite{Polya}.
\begin{Th}
 Let $X$ and $Y$ be i.i.d. centered rv's. Then $X \stackrel{d} = (X+Y)/ \sqrt{2}$  iff
$X$ and $Y$ have the normal distribution with some positive variance.
\end{Th}

The integral test of normality based on this property was proposed by Muliere and Nikitin, see \cite{MN}. Their statistic is asymptotically normal  with the  variance $9\delta^2,$ where
$$
  \delta^2 =\frac{13}{108} -\frac{4}{9\pi}(\arctan\sqrt{\frac{3}{5}} + \frac{1}{2} \arctan{\frac{1}{\sqrt{7}}}) \approx 1.571236\cdot 10^{-3} >0.
$$
The expression for the variance shows the non-trivial character of  the calculations.
The efficiency of this test is very high and equals 0.967 for shift and skew (see \cite{Azz})  alternatives.

We can generalize these findings considering more general characterization, which is the particular case of \cite[Theor.13.7.2]{Kagan}:
\begin{Th}
\label{Gen}
 Let $X$ and $Y$ are centered i.i.d. rv's, and $a$ and $b$ are such constants that $0<a, b<1, a^2 +b^2 =1$ .
 Then $X \stackrel{d}{=} aX +bY$  \,  iff \quad  $X, Y \in N(0,\sigma^2).$
\end{Th}

We can rebuild our statistics using Theorem {\ref{Gen}}, and the result should depend on $a.$  The theory of integral statistic in this generalized setting is developed in \cite{LN}. In particular, the local efficiency of integral test for shift alternative equals to
\begin{multline*}
eff^{*} (a)
    =\left(a-1+\sqrt{1-a^2}\right)^2\Big/\Omega(a), \,  \mbox{where}\\
    \Omega(a) = \Big(\frac{7}{3}\pi-4\arctan{\sqrt{\frac{1+a^2}
    {3-a^2}}}-  4\arctan{\sqrt{\frac{2-a^2}{2+a^2}}} -4\arctan{\sqrt{\frac{1-a^2}{3+a^2}}} - \\ - 4\arctan{\sqrt{\frac{a^2}{4-a^2}}}+4\arctan{\sqrt{\frac{a^2(1-a^2)}{a^4-a^2+4}}}\Big)
    .\nonumber
\end{multline*}

\medskip

The maximum of $eff^{*} (a)$ is 1 but is attained for $a=0$ and $a=1,$ where the test is inconsistent.
 The worst case (quite unexpectedly) is just the Polya case for  $a=\frac{\sqrt{2}}{2}$ with the efficiency $0.966$.  We recommend   $a=\frac{24}{25}$, and $b=\frac{7}{25}$. Then we have $a^2 + b^2 =1,$ and the efficiency is $0.990$, this is a very high value.

 The Kolmogorov type test based on this characterization was studied in \cite{LitNik}. The results are similar but the efficiencies are considerably lower.

 \subsection{Characterization by Shepp property}

In 1964  Shepp \cite{Shepp}  proved that if $X$ and $Y$ are i.i.d., $X, Y \in N(0, \tau^2)$, then the rv
$$k(X,Y):= 2XY/\sqrt{X^2 +Y^2} \in {\mathbb {N}}(0,\tau^2) \quad \text{again}.$$
This statement is usually called the \emph{Shepp property.}

Later Galambos and Simonelli \cite{GS} proved that the
Shepp property characterizes the normal law in some class  $\mathfrak {F_0} $ which consists of such df's  $F$ which satisfy $0< F(0) <1$ and for which $F(x) - F(-x)$ is changing regularly in zero with the exponent 1. They proved the following result
\begin{Th}
 Let $X$ and $Y$ be i.i.d. rv's with common df $F$ from the class  $\mathfrak {F_0}. $ Then the equality in distribution $ X
\stackrel{d}{=} k(X,Y)$ takes place iff  $X \in {\mathbb {N}}(0,\tau^2)$ for some variance $\tau^2 >0.$
\end{Th}

Nikitin and Volkova in \cite{VolNik} constructed tests of normality  based on this characterization and found the efficiencies of corresponding tests.
It turned out that for shift and skew alternatives the efficiencies coincide and are equal in case of integral and supremum tests to
$${ eff_I = \frac{3}{\pi}= 0.955,  \quad  eff_D= \frac{2}{\pi}= 0.637.}$$

\section{Tests of fit for other distributions}

The reader has probably noticed that the majority of characterizations used above for testing exponentiality and normality was formulated in terms of equal distribution of some simple statistics. There arises the question if such characterizations exists for other probability laws and if it is possible to build goodness-of-fit tests based on them. The answer is positive, but the set of corresponding characterizations is more sparse, the calculations are more involved and therefore the whole subject is underdeveloped.

\subsection{Puri-Rubin characterization}

We begin by the characterization of the power law.  We are testing the composite hypothesis
 $$
{  H_0: F \mbox{ is the  df of the power law  so that } F(x)=x^{\mu}, x \in [0,1], \mu > 0,}
 $$
 against general alternatives. We use the characterization which is given in the paper by Puri and Rubin \cite{PR}.
\begin{Th}
\emph{Let $X$ and $Y$ be i.i.d. non-negative rv's with df $F.$
Then the equality $$
X \stackrel{d} = \min(\frac{X}{Y}, \frac{Y}{X})
$$ takes place iff  $X$ and $Y$ have the power distribution.}
\end{Th}

The tests for the power law based on this characterization were build and studied by Nikitin and Volkova \cite{VoNi}.
The efficiencies of integral test are between 0.71 and 0.97, the efficiencies of the Kolmogorov test are between 0.47 and 0.63 depending on the alternative under consideration.

\subsection{Some other laws}

The power law is closely related to the Pareto law, so Obradovic, Jovanovic and Milo\v{s}evic, see \cite{OJM}, were able to use almost the same characterization (by replacing $\min$ by $\max$) when testing for  Pareto law.  Volkova \cite{Voll} introduced and studied some tests of fit for the Pareto distribution based on another characterization.

Goodness-of-fit test for the Cauchy law was build and studied by Litvinova \cite{Litvi}.
She used  the characterization of Ramachandran and Rao \cite{RR}. Its  simplified variant is as follows:
\begin{Th}
{ \it Let $X$ and $Y$ be i.i.d. rv's. Then  $X$ and $\frac13  X  - \frac23 Y $ are identically distributed
 iff  $X$ and $Y$ have the Cauchy df with arbitrary scale factor.}
 \end{Th}

Litvinova in \cite{Litvi} explored the integral test, its local efficiency under the shift alternative turned out to be 0.665.

Some tests of uniformity based on characterizations were developed in \cite{Dud}, \cite{Hash}, \cite{MSMS}. In \cite{Mill} there are interesting efficiency calculations for such tests.

We finish this section by briefly mentioning numerous results on testing goodness-of-fit based on characteristic properties of entropy and Kullback-Leibler information, see \cite{Ari}, \cite{Va}, \cite{Dud}, \cite{Mu}, \cite{Go}, \cite{Nou}, etc. However, there is almost nothing known on efficiencies of new tests, these tests are mainly compared on the basis of simulated power.

\section{Testing of symmetry}

Testing of \emph{symmetry} based on characterizations has been much less explored than goodness-of-fit testing. Consider the classical hypothesis
\begin{equation}
\label{symm}
H_0: 1 - F(x) - F(-x) = 0, \, \, \forall x \in
{\mathbb R}^1,
\end{equation}
against the alternative $H_1$ under which the equality (\ref{symm}) is violated at least in one point. The first step in construction of such tests was done in the crucial paper by Baringhaus and Henze \cite{BH}.

Suppose that $X$ and $Y$ are i.i.d.  rv's with continuous df $F$.  Baringhaus and Henze proved that the equidistribution of rv's $|X|$ and $|\max(X,Y)|$ is valid iff $F$ is symmetric with respect to zero, that is (\ref{symm}) holds. They also proposed suitable Kolmogorov-type  and omega-square type tests of symmetry. Some efficiency calculations for Kolmogorov type test were later performed in \cite{NiBH}, see also \cite{nikiNS}. Integral test of symmetry was next proposed and studied by Litvinova \cite{lit}.

Another characterization of symmetry with respect to 0 belongs to Ahsanullah and was published in \cite{Moe}.
\begin{Th}
{ \textit{Suppose that $X_{1},...,X_{k}, k \geq 2,$ are i.i.d. rv's with absolutely continuous df $F(x).$  Denote $X_{1,k}= \min (X_{1},...,X_{k})$ and $X_{k,k}=\max (X_{1},...,X_{k}).$ Then
$$
|X_{1,k}|  \stackrel{d}= |X_{k,k}|
$$
iff $F$ is symmetric about zero, i.e. $$1-F(x) - F(-x) =0 \quad \forall x \in {\mathbb R}^1.$$ }}
\end{Th}
Subsequently we refer to this result as \emph{ Ahsanullah's characterization of order $k$.}

Nikitin and Ahsanullah  \cite{nah}  published  a paper on  tests of symmetry based on these characterizations
and  their efficiencies. It was found that corresponding tests of symmetry for $ k =2$ and $k = 3$ are asymptotically equivalent to the test of Litvinova and to the Kolmogorov-type test of Baringhaus and Henze. In case of location alternative they are competitive and manifest rather high local Bahadur  efficiency in comparison to many other tests of symmetry. At the same time, higher values of $k, k > 3,$ lead us to different tests with lower values of efficiencies in case of common alternatives. It would be interesting to calculate the efficiencies of such tests for more realistic alternatives, for instance for \emph{skew} alternatives, see \cite{Azz}. First steps in this direction were undertaken in the recent paper \cite{BooNik}.

  Similar research based on a certain modification of Ahsanullah's characterization was done recently by Obradovic and Milo\v{s}evic \cite{MiloO}.
The authors of \cite{MiloO} were able to build corresponding integral test and the test of Kolmogorov type based on their theorem and studied its efficiency.

\section{Directions of further research and perspectives}

\subsection{ Tests based on characterizations of stable laws}

Only three stable laws have explicit densities: normal, Cauchy and L\'evy one-sided density given by the formula

$$
l(x) = \frac{1}{\sqrt{2\pi x^3}} \exp\left(-\frac{1}{2x}\right), \quad x \geq 0.
$$
The tests for normal and Cauchy law based on characterizations were described above.
One of simplest characterizations of the L\'evy law obtained by Ahsanullah and Nevzorov \cite{NAh} looks as follows:\\
\begin{Th}
  Let $X, Y $ and $Z$ be i.i.d. rv's. Then the equality in distribution
$$
X \stackrel{d} = \frac{Y+Z}{4}
$$
takes place  iff $X, Y$ and $Z$ have the one-sided L\'{e}vy distribution  with arbitrary scale factor.
\end{Th}

The tests based on this characterization are unknown. Nothing is known about testing for general stable distributions
using similar characterizations.

 \subsection{ Tests based on characterizations by independence}

The characterization of distributions can be formulated not only in terms of the equidistribution of statistics as in majority of examples given above but also in terms of their independence. Consider as an example the well-known classical result obtained independently  by Kac  \cite{KAC} and Bernstein \cite{BER}  long ago.
\begin{Th}
 If $X$ and $Y$ are independent rv's, then $X+Y$ and $X-Y$  are independent iff $X$ and $Y$ are normal.
 \end{Th}

 As far as we know, this approach is unexplored. Further development of the plot led finally to the famous Skitovich-Darmois theorem \cite{Bryc}, \cite{Kagan} which is also suitable for the construction of tests. We can construct corresponding $U$-empirical distributions and test statistics
 which are more difficult for analysis. Nobody has studied corresponding goodness-of-fit tests.

Another option consists in the well-known result that the independence of $\bar{x}$ and $s^2$ implies normality which was first proved
 by Geary in \cite{Ge}. The same is true for higher central moments. This characterizations were used in a number of papers, see \cite{MK}, \cite{LM}, \cite{Mud},  and a preprint by Thulin \cite{Thu} with power study via simulations. However, there are no calculations of efficiency and analytic comparison with other tests of normality.

First steps in the calculations of efficiency for tests based on characterizations by independence were done recently by Milo\v{s}evic and Obradovic \cite{Milobr}.For instance, they used the following characterization of the exponential law from \cite{Fisz}:
  \begin{Th}
If $X$ and $Y$ are independent i.i.d. random variables with an absolutely continuous distribution and if $\min\{X,Y\}$ and $|X-Y|$ are independent, then both $X$ and $Y$ have exponential distribution with distribution function $F(x)=1-e^{-\lambda x},\;\;x>0,\;\lambda>0$.
\end{Th}
In \cite{Milobr} there are also related results concerning other distributions.

\subsection{Use of empirical integral transforms}

For certain characterizations one can build the test statistics based not on $U$-empirical distributions but on \emph{ empirical transforms, e.g. on empirical characteristic functions or empirical Laplace transforms.}

Let $f_n(t) = n^{-1} \sum_{k=1}^n \exp\{itX_k\}$ be the empirical characteristic function of the sample $X_1, \dots, X_n.$ Then it is clear
that using the Polya characterization ($X \sim (X+Y)/\sqrt{2}$) we have
$$
f_n(t) - f^2_n (t/\sqrt{2})  \approx 0.
$$

Hence the statistics for testing normality of the sample can be
$$
Z_n = \sup_t |  f_n(t) - f^2_n (t/\sqrt{2}) |.
$$
or
$$
W_n = \int_{-\infty}^{\infty}  |  f_n(t) - f^2_n (t/\sqrt{2}) |^2 Q(t) dt,
$$
where $Q$ is some appropriate weight function.

Their asymptotic properties  and efficiencies are unknown. However, the technique of asymptotic analysis of similar statistics was substantially developed in recent years, see, e.g., the papers by Meintanis and Jimenez-Gamero, see \cite{MJ}, \cite{Mei}, \cite{Jim}, \cite{MS}, etc.

The use of empirical Laplace transform with interesting calculation of efficiencies for testing of exponentiality is presented in \cite{MOOO}.

\subsection{Characterizations based on records}

There are many characterizations of distributions based on \emph{record statistics,} see, e.g., \cite{NN}, \cite{De}, \cite{BS}, \cite{SSH},
and many others. Only few of them have been used for construction of goodness-of-fit tests, mainly in the works of Morris and Szynal, see, e.g., \cite{MoSzy}, \cite{MorrS}, \cite{MSzSz}, \cite{Szy} where they essentially used the characterizations based on moments of record values. However, nothing is known about the efficiencies of such tests.

\subsection{Characterizations based on moments}

Some characterizations of distributions are based on their moments or on moments of corresponding order statistics, see, e.g., \cite{Lin}, \cite{LinHu}, \cite{MSz},  \cite{MorMor}, \cite{TooLin}. They can be used for the construction of goodness-of-fit tests, but their efficiencies are unexplored.

\subsection{Multivariate generalizations}

It seems that little or nothing is known about  \emph{multivariate goodness-of-fit tests and multivariate symmetry tests.} One of few exceptions is the recent paper \cite{Milobr}.

\section{Acknowledgements}

This research was supported by grant RFBR No.~16-01-00258 and by grant SPbSU-DFG No.~6.65.37.2017. The author is thankful to the referee
for careful reading of the paper and many useful remarks.

{\footnotesize
 Department of Mathematics and Mechanics, \\
\indent Saint-Petersburg State University, Universitetskaia nab. 7/9,\\
 \indent Saint-Petersburg, 199034, Russia

\medskip

\indent National Research University - Higher School of Economics,\\
\indent Souza Pechatnikov, 16, St.Peters\-burg 190008, Russia

\indent e-mail  \, y.nikitin@spbu.ru   }

\end{document}